\documentclass[reqno]{amsart}

\usepackage{amsmath,amssymb}
\usepackage[psamsfonts]{eucal}

\newtheorem{defi}{Definition}
\newtheorem{teo}{Theorem}

\newtheorem{cor}{Corollary}
\newtheorem{lem}{Lemma}

\title[Full Kostant-Toda lattice]{On the relation between the full Kostant-Toda lattice and multiple orthogonal polynomials}

\author[Barrios]{D. Barrios Rolan\'{\i}a}

\address[Barrios]{Facultad de Inform\'{a}tica,
Universidad Polit\'{e}cnica de Madrid, 28660 Boadilla del Monte, Madrid,
Spain.}
\email[Barrios]{dbarrios@fi.upm.es}

\author[Branquinho]{A. Branquinho}

\address[Branquinho]{CMUC, Department of
Ma\-the\-ma\-tics, University of Coimbra, Largo D. Dinis, 3001-454
Coimbra, Portugal.}
\email[Branquinho]{ajplb@mat.uc.pt}

\author[Foulqui\'e Moreno]{A. Foulqui\'e Moreno}

\address[Foulqui\'e Moreno]{Departamento de
Matem\'atica, Universidade de Aveiro, Campus de Santiago 3810,
Aveiro, Portugal.}
\email[Foulqui\'e Moreno]{foulquie@ua.pt}

\thanks{This research was carried out while on a visit of the first author at University of Coimbra.
The work of this author was supported in part by Direcci\'{o}n General de Investigaci\'{o}n, Ministerio de Educaci\'{o}n y Ciencia, under grant MTM2006-13000-C03-02. The work of the second author was supported by CMUC/FCT. The third author would
like to thank UI Matem\'atica e Aplica\c c\~oes from University of Aveiro.}


\begin{document}

\begin{abstract}
The correspondence between a high-order non symmetric difference operator
with complex coefficients and the evolution of an operator defined by a
Lax pair is established. The solution of the discrete dynamical system is studied, giving explicit
expressions for the resolvent function and, under some conditions, the
representation of the vector of functionals, associated with the solution
for our integrable systems. The method of investigation is based on the evolutions of the matrical
moments.
\end{abstract}

\maketitle

\section{Introduction} \label{sec:1}
We consider the following special \textit{full Kostant-Toda} system,
\begin{equation}
\label{Kostant} \left.
\begin{array}{rcl}
\dot{a}_n & = & b_n-b_{n-1} \\
\dot{b}_n & = & b_n (a_{n+1}-a_n) + c_n - c_{n-1} \\
\dot{c}_n & = & c_n (a_{n+2} - a_n)
\end{array}
\right\}\,,\quad n\in \mathbb{N}\,,
\end{equation}
where the dot means differentiation with respect to $t\in \mathbb{R}$ and we assume $b_0\equiv 0\,,\,c_n\neq 0$.
It is well known that these equations can be written as a \textit{Lax pair} $\dot{J}=\left[J,J_-\right]$, where
$[M,N]=MN-NM$ is the \textit{commutator} of the operators $M$ and $N$, and $J,J_-$ are the operators which matrix
representation is given, respectively,~by
 \begin{equation} \label{J}
J = \left( \begin{matrix}
\phantom{0}  a_1 & 1 & &    \\
 b_1 & a_2 & 1 & \\
c_1& b_2  & a_3  & \ddots \\
0& c_2 & b_3   & \ddots \\
& \ddots &\ddots  &\ddots \\
\end{matrix} \right)\,,\quad
J_- = \left( \begin{matrix}
\phantom{0}  0 & & &   \\
 b_1 & 0 &  &   \\
c_1& b_2  & 0 & \\
0& c_2 & b_3   & \ddots \\
& \ddots &\ddots   & \ddots \\
\end{matrix} \right)\,.
\end{equation}
Here, and in the following, we suppress the explicit $t$-dependence for brevity. Also, we identify an operator
and its matrix representation with respect to the canonical basis. For the sake of simplicity, we only consider a four banded matrix $J$ in this work, but the method can be extended to higher order banded
matrices $J$.

In ~\cite{Geneticas} and ~\cite{*V} the authors considered some special cases of the systems
studied here, and in ~\cite{Ercolani} some finite full Kostant-Toda systems is
considered and solved using bi-orthogonal systems of polynomials.

When $J$ is a bounded operator, then it is possible to define the \textit{resolvent operator},
\begin{equation}
(z I - J)^{-1}=  \sum_{n\geq 0} \frac{J^n  }{z^{n+1}}\,,\quad |z| > \| J \|
 \label{resolvente}
\end{equation}
(see~\cite[Th. 3, p. 211]{Yosida}). We denote by $M_{ij}$ the $2\times 2$ block, of any infinite matrix~$M$,
formed by the entries of rows $2i-1,2i$ and columns $2j-1,2j$. In this way, for each $n\in \mathbb{N}$, $J^n$
can be written as a blocked matrix,
\begin{equation}
J^n = \left(
 \begin{matrix} J^n_{11} & J^n_{12} &  \cdots    \\
J^n_{21} & J^n_{22} & \cdots   \\
\vdots  & \vdots  &   \ddots &
 \end{matrix}
\right)\,. \label{potencia bloques}
 \end{equation}
Moreover, we define $ \mathcal{R}_J(z)$ as the main block of $(z I - J)^{-1}$, this is, $\mathcal{R}_J(z):=(z I
- J)_{11}^{-1}$. Then, from~(\ref{resolvente}) we have
\begin{equation}
\mathcal{R}_J(z)=  \sum_{n\geq 0} \frac{J^n_{11}  }{z^{n+1}}\,,\quad |z| > \| J \|\,.
 \label{resolventefinita}
\end{equation}

As a consequence of the Lax pair representation, for ~(\ref{Kostant}) and other systems, the operator theory is a
useful tool and a remarkable connection between the integrable systems and the approximation theory. Consider
the sequence of polynomials $\{P_n\}$ given by the recurrence relation
\begin{equation}\label{polinomios}
\left.
 \begin{array}{r}
c_{n-1}P_{n-2}(z)+ b_{n}P_{n-1}(z)+(a_{n+1}-z)P_{n}(z)+P_{n+1}(z) =0 \,,\,   n=0,1,\ldots  \\
 P_0(z) = 1 \, ,\quad P_{-1}(z) =P_{-2}(z)=0\,.
\end{array}
\right\}
 \end{equation}
Taking $\mathcal{B}_m(z)= \left( P_{2m}(z), P_{2m+1}(z) \right)^T$, we can rewrite~(\ref{polinomios}) as
\begin{equation}\label{polinomiosvectoriales}
\left.
 \begin{array}{r}
C_{n}\mathcal{B}_{n-1}(z)+ (B_{n+1}-zI_2)\mathcal{B}_{n}(z)+A\mathcal{B}_{n+1}(z) =0\,, \quad n=0,1,\ldots \\
 \mathcal{B}_{-1}(z) = 0 \, ,\quad \mathcal{B}_0(z)=(1,z-a_1)^T
\end{array}
\right\}
 \end{equation}
where
\begin{equation}
A=\left(\begin{array}{cc} 0 & 0\\
1 & 0\end{array}\right)\,,\,
C_n=\left(\begin{array}{cc} c_{2n-1} & b_{2n}\\
0 & c_{2n}
\end{array}\right)\,,\,
B_n=\left(\begin{array}{cc} a_{2n-1} & 1\\
b_{2n-1} & a_{2n}
\end{array}\right)\,,\,n\in \mathbb{N}\,,
\label{matrices} \end{equation} and $C_0$ is an arbitrary $2\times
2$ matrix. In ~(\ref{polinomiosvectoriales}) and in the following,
we denote for the sake of simplicity $(0,0)^T$ and $0\in \mathbb{R}$
in the same way. We recall that the polynomials $P_{n}(z)$ and the
vectorial polynomials $\mathcal{B}_{n}(z)$ depend on $t\in
\mathbb{R}$, when this dependence holds for the coefficients
$a_n,\,b_n,\,c_n$.

Our main goal is to study the solutions of ~(\ref{Kostant}) in
terms of the operator $J$ and its associated vectorial polynomials
$\mathcal{B}_{n}(z)$. We underline that, in a different context
(cf.~\cite[Theorem 2]{primero}), the characterization of solutions
of the integrable system considered there was established in terms
of the derivative of the polynomials associated with $J$. In this
work, using the sequence of vectorial polynomials
$\{\mathcal{B}_{n}\}$, we extend that result. Our first result is
the following.
\begin{teo} \label{equivalencia}
 Assume that the sequence $\{a_n\,,\,b_n\,,\,c_n\}\,,\,n\in
\mathbb{N}\,,$ is uniformly bounded, i.e. there exists $K\in \mathbb{R}_+$ such that
$\max\{|a_n(t)|\,,\,|b_n(t)|\,,\,|c_n(t)|\}\leq M$ for all $n\in \mathbb{N} \text{ and } t\in
\mathbb{R}\,.$ Assume, also, $c_n(t)\neq 0$ for all $n\in \mathbb{N}$ and $t\in \mathbb{R} \, $. Then, the
following conditions are equivalent:
\begin{itemize}
\item[(a)] $\{a_n\,,\,b_n\,,\,c_n\}\,,\,n\in
\mathbb{N}\,,$ is a solution of~(\ref{Kostant}), this is,
\begin{equation}\dot{J}=[J,J_-].\label{Laxpair}
\end{equation}
\item[(b)] For each $n\in \mathbb{N}\cup \{0\} $ we have
\begin{equation}
\frac{d}{dt}J_{11}^n=J_{11}^{n+1}-J_{11}^nB_1+[J_{11}^n, \left(J_-\right)_{11}]. \label{bloques}
\end{equation}
\item[(c)] For all $z\in \mathbb{C}$ such that $|z|>\|J\|$,
 \begin{equation}
\dot{\mathcal{R}}_J(z)=\mathcal{R}_J(z)(zI_2-B_1)-I_2+[\mathcal{R}_J(z),
\left(J_-\right)_{11}] \,,
 \label{derivada bloques}
\end{equation}
where $\mathcal{R}_J(z)$ is given by ~(\ref{resolventefinita}).
\item[(d)]For each $n\in \mathbb{N}\cup \{0\} $, the polynomial $\mathcal{B}_n$ defined by ~(\ref{polinomiosvectoriales}) satisfies
\begin{equation}
\dot{\mathcal{B}}_n(z)=-C_n \mathcal{B}_{n-1}(z)-D_n \mathcal{B}_n(z)\,,\label{derivada polinomios}
\end{equation}
where $D_n=\left(\begin{array}{cc} 0 & 0 \\
b_{2n+1} & 0
\end{array}
\right)$\,.
\end{itemize}
\end{teo}

Let $\mathcal{P}$ be the vector space of polynomials with complex coefficients. It is well known that, given the
recurrence relation~(\ref{polinomios}), there exist two linear moment functionals $u^1,u^2$ from $\mathcal{P}$
to $\mathbb{C}$ verifying, for all $m \in \mathbb{N}$,
\begin{equation}
 \begin{cases}
u^i[z^jP_{2m}]=u^i[z^jP_{2m+1}]=0  \, , \ j=0,1,\ldots ,m-1 \, , \, i=1,2 \\
u^1[z^mP_{2m+1}]=0
 \end{cases}
 \label{funcionales}
\end{equation}
(see~\cite[Th. 3.2]{Iseghem}, see also~\cite{Ismail,Walter}).

We consider the space $\mathcal{P}^2=\{(q_1,q_2)^T: q_i \text{ polynomial},\,i=1,2\}$ and the
space $\mathcal{M}_{2\times 2}$ of $(2\times 2)$-matrices with complex entries.
\begin{defi}
If the functionals $u^1, u^2$ verify ~(\ref{funcionales}), then we say that the function $\mathcal{W}: \mathcal{P}^2 \rightarrow \mathcal{M}_{2\times 2} $ given by
 \begin{equation}
\mathcal{W}
 \left( \begin{matrix} q_1\\q_2 \end{matrix} \right)
 = \left( \begin{matrix}  u^1[q_1]  &
 u^2[q_1] \\
 u^1[q_2]   &  u^2[q_2]
\end{matrix} \right) \, ,
\label{definicion_funcional}
\end{equation}
is a vector of functionals associated with the recurrence relation ~(\ref{polinomios}).
\end{defi}

If $\mathcal{W}$ is a vector of functionals associated with the recurrence relation ~(\ref{polinomios}), then the
following orthogonality relations are verified
\begin{equation}
\mathcal{W}(z^j \mathcal{B}_m) = \mathcal{O}_2\,,\quad j=0,1,\ldots , m-1\,, \label{ortog}
\end{equation}
where $\mathcal{O}_2$ denotes the $2\times 2$ null matrix.

\begin{defi}
A function $\mathcal{W}: \mathcal{P}^2 \rightarrow \mathcal{M}_{2\times 2} $ verifying ~(\ref{ortog}) is called
orthogonality vector of functionals for the recurrence relation ~(\ref{polinomiosvectoriales}).
\end{defi}

Since the above definitions, any vector of functionals associated with the recurrence relation
~(\ref{polinomios}) is always an orthogonality vector of functionals for the recurrence relation
~(\ref{polinomiosvectoriales}). As in the scalar case, it is possible to find more than one orthogonality vector
of functionals. In fact, given a such function $\mathcal{W}: \mathcal{P}^2 \rightarrow \mathcal{M}_{2\times 2}
$, and given any matrix $M\in \mathcal{M}_{2\times 2}$ it is enough to define $\mathcal{W}_M$~as
 \begin{equation}
\mathcal{W}_M \left( \begin{matrix} q_1\\q_2 \end{matrix} \right)
: =\mathcal{W} \left( \begin{matrix} q_1\\q_2 \end{matrix} \right) M \, ,
\label{conmatriz}
 \end{equation}
for having another orthogonality vector of functionals. In the following, we
assume that $\mathcal{W}$ is a fixed vector of functionals associated with the recurrence relation
~(\ref{polinomios}) such that $\mathcal{W} \left( \mathcal{B}_0 \right) $ is an invertible matrix.

We recall that, in ~(\ref{polinomiosvectoriales}), the matrix $C_0$ was arbitrary chosen. In the sequel we assume
 $$
 C_0 = \left( \begin{array}{cc}
1 & 0 \\
-a_1 & 1
\end{array} \right).
 $$
Take $M = \left(
\mathcal{W} \left( \mathcal{B}_0 \right) \right)^{-1} C_0$ and define
$$
\mathcal{U}=\mathcal{W}_M
$$
as in ~(\ref{conmatriz}). Then,
$\mathcal{U}(\mathcal{Q})=\mathcal{W}\left( \mathcal{Q} \right)
\left(\mathcal{W}\left(
\mathcal{B}_0\right)\right)^{-1}C_0$ for any $\mathcal{Q}\in \mathcal{P}^2$ and, in particular,
\begin{equation}
\mathcal{U}(\mathcal{B}_0)=C_0. \label{Udecero} \end{equation}
Moreover, from ~(\ref{polinomiosvectoriales}) and
~(\ref{ortog}),
\begin{equation}
\mathcal{U}(z^m \mathcal{B}_m) =C_m\mathcal{U}(z^{m-1} \mathcal{B}_{m-1})\,,\quad m\in \mathbb{N}
\label{cadena} \end{equation}
(see Lemma \ref{lema*}). Using ~(\ref{Udecero}), ~(\ref{cadena}), and again ~(\ref{ortog}), for each $m\in \mathbb{N}\cup \{0\}$ we
arrive to
\begin{equation}
\mathcal{U}(z^j \mathcal{B}_m) =
\begin{cases}
\mathcal{O}_2 \, , \ j=0,1,\ldots ,m-1 \\
C_mC_{m-1}\cdots C_0 \, , \ j=m\, .
\end{cases}
 \label{U}
\end{equation}

We use the vectors
$ \mathcal{P}_m = \mathcal{P}_m(z) = \left(
z^{2m},z^{2m+1}
 \right)^T$ for each $m\in \mathbb{N}\cup\{0\}$. The following definition extends the corresponding to the scalar case.
\begin{defi}
\label{definicion momentos} For each $m=0,1,\ldots $, the matrix $\mathcal{U}
\left( z^m \mathcal{P}_0 \right)$ is
called moment of order $m$ for the vector of functionals $\mathcal{U}$.
\end{defi}

In particular, since $\mathcal{B}_0=C_0\mathcal{P}_0$, we have
\begin{equation}
\mathcal{U} \left( \mathcal{P}_0 \right) = I_2 \label{momento0}
\end{equation}
(see ~(\ref{Udecero})).

We know $\mathcal{U}=\mathcal{U}\{t\}$ depends on $t$, besides this dependence is not explicitly written (as we said above). Then, it is possible to define the derivative of $\mathcal{U}$ as usual,
$$
\frac{d\mathcal{U}}{dt}:\mathcal{P}^2 \longrightarrow \mathcal{M}_{2\times 2}
$$
such that, for each $\mathcal{B}\in\mathcal{P}^2$,
$$
\frac{d\mathcal{U}}{dt}(\mathcal{B})=\lim_{\Delta t\to 0}\frac{\mathcal{U}\{t+\Delta t\}(\mathcal{B})-
\mathcal{U}\{t\}(\mathcal{B})}{\Delta t}\,.
$$
Obviously, the usual properties for this kind of operators are verified. In particular,
\begin{equation}
\frac{d}{dt} \left( \mathcal{U} (\mathcal{B}) \right)
= \frac{d\mathcal{U}}{dt}(\mathcal{B}) +
\mathcal{U} ( \dot{\mathcal{B}} )
 \,, \quad \forall \mathcal{B}\in\mathcal{P}^2 \,.
\label{*}
\end{equation}
We use $\frac{d\mathcal{U}}{dt}$ and ~(\ref{*}) below. Also, we will use the matrix function $\mathcal{R}_J$ given in ~(\ref{resolventefinita}). We define the \textit{generating function of the moments} as
\begin{equation}
\mathcal{F}_J(z)=C_0^{-1}\mathcal{R}_J(z)C_0\,,\quad |z|>\|J\|\,.
\label{generadora}
\end{equation}

Next, we have our second main result, related with Theorem~\ref{equivalencia}. More precisely speaking, we will
see that Theorem~\ref{equivalencia} follows directly from Theorem~\ref{teorema2}.
\begin{teo} \label{teorema2}
In the conditions of Theorem~\ref{equivalencia}, assume $\dot{a}_1=b_1$. Let $\mathcal{U}$ be given by~(\ref{U}). Then, the following assertions are
equivalent:
\begin{itemize}
\item[(e)] $\{a_n\,,\,b_n\,,\,c_n\}\,,\,n\in
\mathbb{N}\,,$ is a solution of~(\ref{Kostant}), this is, ~(\ref{Laxpair}) holds.
\item[(f)] For each $n=0,1,\ldots , $ we have
\begin{equation}
\frac{d}{dt}\mathcal{U} \left( z^n\mathcal{P}_0 \right)
 =
 \mathcal{U}\left(z^{n+1}\mathcal{P}_0\right)-
\mathcal{U}\left(z^n\mathcal{P}_0\right)\mathcal{U}\left(z\mathcal{P}_0\right)\,.\label{derivada momentos}
\end{equation}
\item[(g)] For all $\zeta\in \mathbb{C}$ such that $|\zeta|>\|J\|$,
\begin{equation}
\dot{\mathcal{F}}_J(\zeta) =\mathcal{F}_J(\zeta)\left(
\zeta I_2-\mathcal{U}(z\mathcal{P}_0)\right)-I_2\,,\label{derivada generatriz}
\end{equation}
being $\mathcal{F}_J$ the generating function defined in ~(\ref{generadora}).
\item[(h)] For all $\mathcal{B}\in \mathcal{P}^2$ we have
\begin{equation}
\left(\frac{d}{dt}\mathcal{U}\right)(\mathcal{B})=\mathcal{U}(z\mathcal{B})-\mathcal{U}(\mathcal{B})
\mathcal{U}(z\mathcal{P}_0)\,. \label{derivada funcional}
\end{equation}
\item[(i)]For each $n=0,1,\ldots ,$ we have ~(\ref{derivada polinomios})\,.
\end{itemize}
\end{teo}

Moreover, we have other consequences of Theorem~\ref{equivalencia} and Theorem~\ref{teorema2}.

In the next result, $\int_0^tf(s)ds$ is understood to be the solution $X=X(t)$ of the Cauchy problem
$$\left.\begin{array}{r}
\dot{X}  = f(t)\\
X(0) = 0
\end{array}\right\}
$$
in the suitable space. It is well-known that, in our conditions, there exists a unique solution of this problem (see, for instance, ~\cite{V} and ~\cite{Vilenkin}).
\begin{cor}\label{coro1}
Under the conditions of Theorem~\ref{equivalencia}, if
$\{a_n\,,\,b_n\,,\,c_n\}\,,\,n\in
\mathbb{N}\,,$ is a solution of~(\ref{Kostant}), then
\begin{equation}
\mathcal{R}_J(z)= \exp \left( zt \right)
C_0M(t,z) \left( N(t) \right)^{-1}\,,
\label{representacion}
\end{equation}
where
$$
N(t) = \left( \begin {array}{cc}
\exp \left( \int_0^ta_1ds \right) &
\exp \left( \int_0^ta_1ds \right)
\int_0^t \exp\left(-\int_0^s(a_2-a_1)dr\right) ds \\
0 & \exp \left( \int_0^ta_2ds \right)
\end{array}
\right)\,,
$$
$$
M(t,z) = -\int_0^t\exp \left( -zs \right) C_0^{-1}N(s)ds + \left( C_0(0) \right)^{-1}
\mathcal{R}_0(z)
$$
(here, $C_0(0)$ and $\mathcal{R}_0(z)$ are, respectively, $C_0$ and $\mathcal{R}(z)$ for $t=0$).
\end{cor}

Given a linear functional $u:\mathcal{P}\rightarrow \mathbb{C}$, we may define the new functional  $e^{zt}u:\mathcal{P}\rightarrow \mathbb{C}$ as
\begin{equation}
\left( e^{zt} u \right) [z^j] = \sum_{k\geq 0}\frac{t^k}{k!}u[z^{k+j}]\,.
\label{exponencial}
\end{equation}
We denote by $\mathcal{U}_0=(u_0^1,u_0^2)^T$ our vector of functionals $\mathcal{U}$ when $t=0$ and, similarly, by $J_0$ the triangular matrix given in ~(\ref{J}) when $t=0$. If $J_0$ is a bounded matrix, then
$$
|u^i_0[z^{k+j}]|\leq m_{ij}\|J_0\|^{k+j}\,,\quad i=1,2\,,\quad k,j=0,1,\ldots ,
$$
and the right-hand side of ~(\ref{exponencial}) is well-defined for $u=u^i_0\,,\,i=1,2\,$ (see \cite[Theorem~4]{Kaliaguine}). In this case, we can define the vector of functionals $e^{zt} \, \mathcal{U}_0$ as
$$
\left( e^{zt} \, \mathcal{U}_0 \right)
 \left( \mathcal{B} \right)
 = \left( \begin{matrix} \left( e^{zt} u_0^1 \right) [q_1]  &
\left( e^{zt}u_0^2 \right) [q_1] \\
 \left(e^{zt}u_0^1\right)[q_2] & \left( e^{zt}u_0^2 \right) [q_2]
\end{matrix} \right)\,
$$
for each $\mathcal{B}=(p,q)^T\in\mathcal{P}^2$. As in \cite[Theorem 3]{Dynamics}, we give a possible representation for the vector of functionals $\mathcal{U}$.
\begin{cor}
In the conditions of Theorem~\ref{teorema2}, and with the notation used in~(\ref{conmatriz}), assume that the vector of functionals $\mathcal{U}$ verifies
\begin{equation}
\mathcal{U}=
\left( e^{zt} \, \mathcal{U}_0 \right)_M
\label{Uexponencial}
\end{equation}
for some $M\in \mathcal{M}_{2\times 2}$. Then,
$\{a_n\,,\,b_n\,,\,c_n\}\,,\,n\in
\mathbb{N}\,,$ is a solution of~(\ref{Kostant}).
\label{corolario2}
\end{cor}

In section~\ref{sec:2} we show that the study of the system~\eqref{Kostant} can
be reduced to the evolution of the main block of the matrix $J^n $ i.e. $ J_{11}^n$.
In section~\ref{sec:3} we prove Theorem~\ref{teorema2}. The main idea is to express~\eqref{Kostant} in
terms of the evolution of the moments.
In section~\ref{sec:4} we prove Theorem~\ref{equivalencia}. The main feature of this result is the
connections between the resolvent function, $\mathcal{R}_J$, and the
evolution of the systems of vector polynomial, $\{ \mathcal{B}_n \}$.

\section{Auxiliary results} \label{sec:2}

Next lemma was used for proving ~(\ref{cadena}).

\begin{lem}\label{lema*}
Let $\mathcal{W}$ be a vector of functionals  associated with ~(\ref{polinomios}). Then
\begin{equation}
\mathcal{W}(A_1\mathcal{Q}_1 + A_2\mathcal{Q}_2) = A_1\mathcal{W}(\mathcal{Q}_1 ) +
A_2\mathcal{W}(\mathcal{Q}_2)\label{propiedad}
\end{equation}
is verified for any $\mathcal{Q}_1,\mathcal{Q}_2\in \mathcal{P}^2$ and $  A_1,A_2 \in \mathcal{M}_{2\times 2}
\,$.
\end{lem}

\textbf{Proof}.- It is sufficient to take into account that $\mathcal{W}: \mathcal{P}^2 \rightarrow \mathcal{M}_{2\times 2} $
is given by
$$
\mathcal{W}
 \left( \begin{matrix} q_1\\q_2 \end{matrix} \right)
 = \left(\begin{matrix} u[q_1] &
 v[q_1] \\
 u[q_2] & v[q_2]
 \end{matrix} \right)\,$$
when $u,v: \mathcal{P} \rightarrow \mathbb{C}$ are linear functionals. \hfill $\square$

\begin{lem}\label{lema**}
The orthogonality vector of functionals
$\mathcal{U}:\mathcal{P}^2 \longrightarrow \mathcal{M}_{2\times 2}$ is determined by~(\ref{U}). This is, $\mathcal{U}$ is the unique vector of functionals associated with the recurrence relation ~(\ref{polinomios}) verifying~(\ref{U}).
\end{lem}

\textbf{Proof}.-
Given $\left(q_1,q_2\right)^T\in \mathcal{P}^2$, for each $i=1,2$ we can write
$$
q_i(z) = \sum_{k=1}^2\alpha^0_{ik} P_{k-1}(z) + \sum_{k=1}^2 \alpha^1_{ik} P_{1+k}(z) + \cdots + \sum_{k=1}^2
\alpha^m_{ik} P_{2m+k-1}(z) \, , \ \alpha^j_{ik} \in \mathbb{C} \, ,
 $$
where $m=\max\{m_1,m_2\}$ and $\deg (q_i)\leq 2m_i+1$ (we understand $\alpha^j_{ik}=0$ when $j>m_i$). This is,
\begin{equation}
\left(q_1,q_2\right)^T=\sum_{j=0}^mE_j \, \mathcal{B}_j\,, \label{combinacion} \end{equation} being
$E_j=\left(\alpha^j_{ik}\right)\in \mathcal{M}_{2\times 2}\,,\,j=0,\ldots ,m$. From~(\ref{propiedad}), if
$\widetilde{\mathcal{U}}:\mathcal{P}^2 \longrightarrow \mathcal{M}_{2\times 2}$ is a vector of functionals associated with the recurrence relation ~(\ref{polinomios}), then
\begin{equation}
\widetilde{\mathcal{U}}\left(\begin{array}{c}q_1\\q_2\end{array}\right)=\sum_{j=0}^mE_j \,\widetilde{\mathcal{U}} \left(\mathcal{B}_j\right)\,.
\label{otra}
\end{equation}
Moreover, if $\widetilde{\mathcal{U}}$ verifies~(\ref{U}) we have
$
\widetilde{\mathcal{U}} \left(\mathcal{B}_j\right)=\mathcal{U} \left(\mathcal{B}_j\right)=0\,,\,j=1,2,\ldots \,,
$
and also
$
\widetilde{\mathcal{U}} \left(\mathcal{B}_0\right)=\mathcal{U} \left(\mathcal{B}_0\right)=C_0\,.
$
Therefore, from ~(\ref{otra}) we arrive to
$\widetilde{\mathcal{U}}=\mathcal{U}\,.
$
\hfill $\square$

Next result shows that it is possible to recover the entries of matrix $J$ using the orthogonality vector of functionals $\mathcal{U}$.

\begin{lem}
\label{lema***}
The entries of the matrix $J$ are determined by the sequence $\{\mathcal{B}_n\}$ of vectorial polynomials.
\end{lem}

\textbf{Proof}.- The entries of $J$ are determined by the blocks $C_n\,,\,B_n\,,\,n\in \mathbb{N}$. Then, it is sufficient to express these blocks in terms of $\{\mathcal{B}_n\}$. Since $C_k\,,\,k=0,1,\ldots \,,$ are invertible matrices, also $\mathcal{U}\left(z^k \mathcal{B}_k \right)$ is invertible and, from ~(\ref{U}),
$$
C_n=\mathcal{U}\left(z^n\mathcal{B}_n \right)\left(\mathcal{U}\left(z^{n-1} \mathcal{B}_{n-1} \right)\right)^{-1}\,,\quad n=1,2,\ldots \,.
$$
On the other hand, from ~(\ref{polinomiosvectoriales}) and ~(\ref{U}) we deduce
$$
C_n\mathcal{U}\left(z^n \mathcal{B}_{n-1} \right)+B_{n+1}\mathcal{U}\left(z^n \mathcal{B}_n \right)-\mathcal{U}\left(z^{n+1} \mathcal{B}_n \right)=0\,,\quad n=0,1,\ldots\,.
$$
Then, for $n\in \mathbb{N}$ we have
$$
B_{n}=\left(\mathcal{U}\left(z^{n} \mathcal{B}_{n-1} \right)-C_{n-1}\mathcal{U}\left(z^{n-1} \mathcal{B}_{n-2} \right)\right)\left(\mathcal{U}\left(z^{n-1} \mathcal{B}_{n-1} \right)\right)^{-1} \, ,
$$
and the result follows. \hfill $\square$

Next, we determine the expression of the moment $\mathcal{U}\left(\mathcal{P}_n\right)=\mathcal{U}\left(x^n \,
\mathcal{P}_0\right)$ in terms of the matrix~$J$.
\begin{lem}\label{expresion momentos}
For each $n=0,1,\ldots $ we have \begin{equation}
 \displaystyle
 \mathcal{U} \left(z^n \, \mathcal{P}_0\right)=
C_0^{-1} J^{n}_{11} C_0\,. \label{momentos}
\end{equation}
 \end{lem}
\textbf{Proof}.- We know that
$\mathcal{U}\left(\mathcal{P}_0\right)=I_2$
(see~(\ref{momento0})), then the result is verified for $n=0$.

Another way to express ~(\ref{polinomiosvectoriales}) is
\begin{equation*}
J
\left(  \begin{matrix} \mathcal{B}_0(z) \\
\mathcal{B}_1 (z)\\
\vdots \\
\end{matrix}\right) =z \left(  \begin{matrix} \mathcal{B}_0 (z)\\
\mathcal{B}_1(z) \\
\vdots
\end{matrix}\right)\,.
 \end{equation*}
Thus,
\begin{equation} \label{potencias}
 J^n
\left(  \begin{matrix} \mathcal{B}_0(z) \\
\mathcal{B}_1 (z)\\
\vdots \\
\end{matrix}\right) = z^n \left(  \begin{matrix} \mathcal{B}_0 (z)\\
\mathcal{B}_1 (z)\\
\vdots
\end{matrix}\right)\,,\quad n\in \mathbb{N}\,.
\end{equation}
Comparing the first rows in ~(\ref{potencias}), and taking into account ~(\ref{potencia bloques}) and the form of
$J$,
$$\displaystyle \sum_{i\geq 1}J^n_{1i} \, \mathcal{B}_{i-1}(z) =
J^n_{11}\mathcal{B}_{0}(z)+J^n_{12}\mathcal{B}_1(z)=z^n \, \mathcal{B}_0(z)\,.
 $$
Then, from ~(\ref{U}),
$$
\mathcal{U}\left(z^n \, \mathcal{B}_0\right)= J^n_{11}\mathcal{U}\left(\mathcal{B}_{0}\right)\,,
$$
this is,
$$
C_0\,\mathcal{U}\left(z^n\mathcal{P}_{0}\right)=J^n_{11}C_0
$$
(see ~(\ref{Udecero}) and ~(\ref{momento0})), which is ~(\ref{momentos}). \hfill $\square$

The following result concerns to solutions associated with the matrix $J$, non necessarily bounded.

\begin{lem}
If the sequence $\{a_n\,,\,b_n\,,\,c_n\}\,,\,n\in \mathbb{N}\,,$ is a solution of~(\ref{Kostant}),
then ~(\ref{bloques}) is verified. \label{lema2}
\end{lem}

\textbf{Proof}.- Under the given conditions, we know that ~(\ref{Laxpair}) holds. Then, it is very easy to verify
$$
\frac{d}{dt}J^n=J^nJ_--J_-J^n
$$
and, with the established notation,
$$
\frac{d}{dt}J^n_{11}=\left(J^nJ_-\right)_{11}-\left(J_-J^n\right)_{11}
$$
 From ~(\ref{J}) and ~(\ref{potencia bloques}),
$$
\left.\begin{array}{rcl} \left(J^nJ_-\right)_{11} & = & J^n_{11}\left(J_-\right)_{11}+
J^n_{12}\left(J_-\right)_{21}\\
\left(J_-J^n\right)_{11} & = & \left(J_-\right)_{11}J^n_{11}\,.
\end{array}
\right\}
$$
Then, \begin{equation} \dot{J}^n_{11}=J^n_{11}\left(J_-\right)_{11}- \left(J_-\right)_{11}J^n_{11} +
J^n_{12}\left(J_-\right)_{21}\,. \label{4} \end{equation} On the other hand, $
J^{n+1}_{11}=J^n_{11}J_{11}+J^n_{12}\left(J_-\right)_{21}\,, $ which, joint with ~(\ref{4}), \linebreak goes
to~(\ref{bloques}). \hfill $\square$

\section{Proof of Theorem~\ref{teorema2}} \label{sec:3}

In the first place, we show (e)$\Rightarrow$ (f). Assume that ~(\ref{Laxpair}) holds. Since Lemma \ref{lema2}, we
have ~(\ref{bloques}). Moreover, due to Lemma \ref{expresion momentos},
$J^n_{11}=C_0\mathcal{U}(z^n\mathcal{P}_0)C_0^{-1}\,,\,n\in \mathbb{N}\,,$ and, in particular,
$B_1=J_{11}=C_0\mathcal{U}(z\mathcal{P}_0)C_0^{-1}$. Also,
$C_0^{-1}\left(J_-\right)_{11}=\left(J_-\right)_{11}C_0=\left(J_-\right)_{11}$. Then, from ~(\ref{4}),
\begin{multline}
\frac{d}{dt}\left(C_0\mathcal{U}(z^n\mathcal{P}_0)C_0^{-1}\right) =
C_0\mathcal{U}(z^{n+1}\mathcal{P}_0)C_0^{-1}-C_0\mathcal{U}(z^n\mathcal{P}_0)\mathcal{U}(z\mathcal{P}_0)  C_0^{-1}  \\
 +  C_0\mathcal{U}(z^n\mathcal{P}_0)\left(J_-\right)_{11}
-\left(J_-\right)_{11}\mathcal{U}(z^n\mathcal{P}_0)C_0^{-1} \label{5}
\end{multline}
On the other hand, taking derivatives (and denoting by $\dot{C}_0^{-1}$ the derivative of $C_0^{-1}$),
\begin{multline}
\frac{d}{dt}\left(C_0\mathcal{U}(z^n\mathcal{P}_0)C_0^{-1}\right)  =
\dot{C}_0\mathcal{U}(z^n\mathcal{P}_0)C_0^{-1}+C_0\mathcal{U}(z^n\mathcal{P}_0)\dot{C}_0^{-1}  \\
  +  C_0\frac{d}{dt}\left(\mathcal{U}(z^n\mathcal{P}_0)\right)C_0^{-1} \phantom{olao} \label{6}
\end{multline}
Since $\dot{a}_1=b_1$, we can verify $\dot{C}_0=-\left(J_-\right)_{11}\,,\,
\dot{C}_0^{-1}=\left(J_-\right)_{11}$. Hence, comparing the right hand sides of ~(\ref{5}) and ~(\ref{6}) we
arrive to ~(\ref{derivada momentos}) (we recall that $C_0$ is an invertible matrix).

Now, we prove (f)$\Rightarrow$ (g). Since ~(\ref{resolventefinita}) and Lemma~\ref{expresion momentos},
\begin{equation}
\mathcal{F}_J(\zeta)=\sum_{n\geq 0}\frac{\mathcal{U}(z^n\mathcal{P}_0)}{\zeta^{n+1}}\,,\quad |\zeta|>\|J\|\,.
\label{F}
\end{equation}
Then, from ~(\ref{derivada momentos}),
\begin{eqnarray*}
\dot{\mathcal{F}}_J(\zeta) & = & \sum_{n\geq 0}\frac{\mathcal{U}(z^{n+1}\mathcal{P}_0)-
\mathcal{U}(z^n\mathcal{P}_0)\mathcal{U}(z\mathcal{P}_0)}{\zeta^{n+1}} \\
& = & \zeta\sum_{n\geq 0}\frac{\mathcal{U}(z^{n+1}\mathcal{P}_0)}{\zeta^{n+2}}-\sum_{n\geq
0}\frac{\mathcal{U}(z^n\mathcal{P}_0)}{\zeta^{n+1}}\mathcal{U}(z\mathcal{P}_0) \\
& = & \zeta\left( \mathcal{F}_J(\zeta)-\frac{1}{\zeta}I_2
\right)-\mathcal{F}_J(\zeta)\mathcal{U}(z\mathcal{P}_0)\,.
\end{eqnarray*}
This is, ~(\ref{derivada generatriz}) is verify.

Given ~(\ref{derivada generatriz}), we are going to obtain the
derivative of the vector of functionals $\mathcal{U}$. For doing
this, we use the linearity of $\mathcal{U}$ and the convergence of
the series,
\begin{equation}
\mathcal{F}_J(\zeta)=\mathcal{U}\left(\sum_{n\geq 0}\frac{z^n}{\zeta^{n+1}}\mathcal{P}_0\right)
=\mathcal{U}\left(\frac{1}{\zeta-z}\mathcal{P}_0\right)\,,\quad |\zeta|>\|J\|\,.
\label{3}
\end{equation}
(Here and in the next expressions, as usual, $\mathcal{U}=\mathcal{U}_z$ is the vector of functionals defined on the closure of the space $\mathcal{P}^2$ of vectorial polynomials $(q_1,q_2)^T$ in the variable $z$.)

 From ~(\ref{3}) and ~(\ref{derivada generatriz}),
\begin{eqnarray}
\frac{d}{dt}\mathcal{U}\left(\frac{1}{\zeta-z}\mathcal{P}_0\right)& = &
\mathcal{U}\left(\frac{1}{\zeta-z}\mathcal{P}_0\right)\left(\zeta I_2-\mathcal{U}\left(
z\mathcal{P}_0\right)\right) -I_2\nonumber \\
& = & \mathcal{U}\left(\left(1+\frac{z}{\zeta-z}\right)\mathcal{P}_0\right)-
\mathcal{U}\left(\frac{1}{\zeta-z}\mathcal{P}_0\right)\mathcal{U}\left(z\mathcal{P}_0\right)-
I_2 \nonumber \\
& = & \mathcal{U}\left(\frac{z}{\zeta-z}\mathcal{P}_0\right)-
\mathcal{U}\left(\frac{1}{\zeta-z}\mathcal{P}_0\right)\mathcal{U}\left(z\mathcal{P}_0\right)\,.
\label{VIII}
\end{eqnarray}
Define the vectors of functionals $\mathcal{U}_1\,,\,\mathcal{U}_2:\mathcal{P}^2\longrightarrow
\mathcal{M}_{2\times 2}$ such that
\begin{equation}
\left.\begin{array}{rcl} \mathcal{U}_1(\mathcal{B}) & = & \mathcal{U}(z\mathcal{B})\\
\mathcal{U}_2(\mathcal{B}) & = &\mathcal{U}(\mathcal{B}) \mathcal{U}(z\mathcal{P}_0)
\end{array}
\right\} \label{nuevos}
\end{equation}
for each $\mathcal{B}\in \mathcal{P}^2$. We remark that
$\frac{1}{\zeta-z}\mathcal{P}_0$ do not depend on $t\in
\mathbb{R}$. In ~(\ref{VIII}), denoting
$\dot{\mathcal{U}}=\frac{d}{dt}\mathcal{U}$, we have
$\dot{\mathcal{U}}=\mathcal{U}_1-\mathcal{U}_2$ over
$\frac{1}{\zeta-z}\mathcal{P}_0$, being
$$
\mathcal{U}\left(\frac{1}{\zeta-z}\mathcal{P}_0\right)=\frac{1}{\zeta}\mathcal{U}\left(\mathcal{P}_0\right)
+\frac{1}{\zeta^2}\mathcal{U}\left(z\mathcal{P}_0\right)+\cdots\,,\quad |\zeta|>\|J\|\,.
$$
Hence, we have $\dot{\mathcal{U}}=\mathcal{U}_1-\mathcal{U}_2$ over $\mathcal{P}^2$, this is, we have
~(\ref{derivada funcional}).

For proving (h)$\Rightarrow$ (i), as in ~(\ref{combinacion}), $\dot{\mathcal{B}}_n$ can be written in terms of
the sequence~$\{\mathcal{B}_n\}$,
\begin{equation}
\dot{\mathcal{B}}_n(z)=D^{(n)}_0\mathcal{B}_0(z)+D^{(n)}_1\mathcal{B}_1(z)+\cdots +D^{(n)}_n\mathcal{B}_n(z)
\label{IX}
\end{equation}
If $n=0,1$, the above expression is
\begin{equation}
\dot{\mathcal{B}}_n(z)=D^{(n)}_{n-1}\mathcal{B}_{n-1}(z)+D^{(n)}_n\mathcal{B}_n(z) \label{X}
\end{equation}
Let $n\geq 2$ be fixed. We are going to show that ~(\ref{X}) holds, also, for $n$. Due to the orthogonality, from
~(\ref{IX}),
$$
\mathcal{U}(\dot{\mathcal{B}}_n)=D^{(n)}_{0}\mathcal{U}\left(\mathcal{B}_{0}\right)\,.
$$
In fact, using ~(\ref{derivada funcional}),
\begin{eqnarray*}
\mathcal{O}_{2}=\frac{d}{dt}\left(\mathcal{U}({\mathcal{B}}_n)\right)& = & \dot{\mathcal{U}}({\mathcal{B}}_n)+
\mathcal{U}(\dot{{\mathcal{B}}}_n)\\
& = &\mathcal{U}(z{\mathcal{B}}_n)-\mathcal{U}(\mathcal{B}_n)\mathcal{U}
(z{\mathcal{P}}_0)+D^{(n)}_{0}\mathcal{U}\left(\mathcal{B}_{0}\right)\,.
\end{eqnarray*}
Thus, $D^{(n)}_{0}=\mathcal{O}_{2}$. We proceed by induction on $n$, assuming
$$
D^{(n)}_{0}=\cdots =D^{(n)}_{j-1}=\mathcal{O}_2
$$
for a fixed $j<n-1$. Using ~(\ref{IX}) and, again, ~(\ref{derivada funcional}) and ~(\ref{U}),
\begin{eqnarray*}
\mathcal{O}_{2}=\frac{d}{dt}\mathcal{U}(z^j\mathcal{B}_n)& = &
\left(\frac{d}{dt}\mathcal{U}\right)(z^j\mathcal{B}_n)+
\mathcal{U}(z^j\dot{{\mathcal{B}}}_n)\\
& = &\mathcal{U}(z^{j+1}{\mathcal{B}}_n)-\mathcal{U}(z^j\mathcal{B}_n)\mathcal{U}
(z{\mathcal{P}}_0)+D^{(n)}_{j}\mathcal{U}\left(z^j\mathcal{B}_{j}\right)\\
& = & D^{(n)}_{j}\mathcal{U}\left(z^j\mathcal{B}_{j}\right)\,,
\end{eqnarray*}
where $\mathcal{U}\left(z^j\mathcal{B}_{j}\right)$ is an invertible matrix. Thus, $D^{(n)}_{j}=\mathcal{O}_{2}$
and ~(\ref{X}) is verified for any $n\in \mathbb{N}$.

Our next purpose is to determine $D^{(n)}_{j}\,,\,j=n-1,n$. From ~(\ref{X}),
$$
\mathcal{U}(z^{n-1}\dot{\mathcal{B}}_{n})=D^{(n)}_{n-1}\mathcal{U}\left(z^{n-1}\mathcal{B}_{n-1}\right)\,.
$$
Then, because of ~(\ref{derivada funcional}) and ~(\ref{U}),
$$
\mathcal{O}_{2}=\frac{d}{dt}\mathcal{U}\left(z^{n-1}\mathcal{B}_{n}\right)=\mathcal{U}\left(z^{n}\mathcal{B}_{n}\right)-
\mathcal{U}\left(z^{n-1}\mathcal{B}_{n}\right)\mathcal{U}\left(z\mathcal{P}_0\right)+ D^{(n)}_{n-1}
\mathcal{U}\left(z^{n-1}\mathcal{B}_{n-1}\right)
$$
and, therefore,
$$
D^{(n)}_{n-1}=-C_n\,.
$$

On the other hand, writing
\begin{equation}
\mathcal{B}_n(z)=\sum_{i=0}^n F^{(n)}_i\mathcal{P}_i(z) \label{XI}
\end{equation}
and comparing the coefficient of $z^{2n}$ and $z^{2n+1}$ in both sides of ~(\ref{XI}), we obtain
$$
F^{(n)}_n=\left(\begin{array}{cc} 1 & 0\\ f_n & 1
\end{array}
\right)\,.
$$
Moreover, taking derivatives in ~(\ref{XI}) and comparing with ~(\ref{X}), we see $D^{(n)}_{n}=\dot{F}^{(n)}_n$
or, what is the same,
$$
D^{(n)}_n=\left(\begin{array}{cc} 0 & 0\\ d_n & 0
\end{array}
\right)\,,
$$
where we need to determine $d_n$. From ~(\ref{derivada funcional}) and ~(\ref{X}),
 \begin{multline*}
D^{(n)}_n\mathcal{U}\left(z^{n}\mathcal{B}_{n}\right) \\ =
\left(\frac{d}{dt}\mathcal{U}\left(z^{n}\mathcal{B}_{n}\right)-\mathcal{U}\left(z^{n+1}\mathcal{B}_{n}\right)
+\mathcal{U}\left(z^{n}\mathcal{B}_{n}\right)\mathcal{U}\left(z\mathcal{P}_{0}\right)\right)+
C_n\mathcal{U}\left(z^{n}\mathcal{B}_{n-1}\right) .
 \end{multline*}
Then, using ~(\ref{U}) and ~(\ref{polinomiosvectoriales}),
\begin{multline*}
D^{(n)}_n  =  \left(\frac{d}{dt}\left(C_nC_{n-1}\cdots C_0\right)\right)\left(C_nC_{n-1}\cdots
C_0\right)^{-1} \\
  -B_{n+1}+ \left(C_nC_{n-1}\cdots C_0\right) \mathcal{U}\left(z\mathcal{P}_{0}\right)\left(C_nC_{n-1}\cdots
C_0\right)^{-1} \, ,
\end{multline*}
thus
\begin{multline}
D^{(n)}_n +B_{n+1}  =  \left(\frac{d}{dt}\left(C_nC_{n-1}\cdots C_1\right)\right)\left(C_nC_{n-1}\cdots C_1\right)^{-1}
   \\
 +\left(C_nC_{n-1}\cdots C_1\right) \left( \dot{C}_0C_0^{-1} +J_{11}\right)
\left(C_nC_{n-1}\cdots C_1\right)^{-1} \, .
 \label{triangular}
\end{multline}
The matrix $C_nC_{n-1}\cdots C_1$ is upper triangular. Moreover, because of $\dot{a}_1=b_1$ also
$\dot{C}_0C_0^{-1} +J_{11}$ is an upper triangular matrix and, then, the matrix in the left hand side of
~(\ref{triangular}) is upper triangular and, consequently, $d_n=b_{2n+1}$.

Finally, we show (i) $\Rightarrow$ (e). Taking derivatives in ~(\ref{polinomiosvectoriales}),
\begin{multline*}
\dot{C}_{n}\mathcal{B}_{n-1}(z) + \dot{ B}_{n+1}\mathcal{B}_{n}(z)   \\
 + C_{n}\dot{\mathcal{B}}_{n-1}(z)  +  (B_{n+1}-zI_2)\dot{\mathcal{B}}_{n}(z)+A\dot{\mathcal{B}}_{n+1}(z) = 0
\,,\quad n= 0,1, \ldots
\end{multline*}
Using ~(\ref{derivada polinomios}) and taking into account $AD_{n+1}=D_nA=\mathcal{O}_{2}$,
\begin{multline}
\dot{C}_{n}\mathcal{B}_{n-1}(z) + \dot{ B}_{n+1}\mathcal{B}_{n}(z) + C_{n}
\left( -C_{n-1}\mathcal{B}_{n-2}(z)+D_{n-1}\mathcal{B}_{n-1} (z) \right)  \\
 + (B_{n+1}-zI_2) \left( -C_{n}\mathcal{B}_{n-1}(z)+D_n\mathcal{B}_n(z) \right)
 - A C_{n+1}\mathcal{B}_{n}(z) = 0 \, .
  \label{8}
\end{multline}
Using, again, ~(\ref{polinomiosvectoriales}) for eliminating the explicit expression in $z$,
$$
\left.\begin{array}{rcl} z\mathcal{B}_n(z) & = &C_{n}\mathcal{B}_{n-1}(z)+
B_{n+1}\mathcal{B}_{n}(z)+A\mathcal{B}_{n+1}(z)\\
z\mathcal{B}_{n-1}(z) & = &C_{n-1}\mathcal{B}_{n-2}(z)+ B_{n}\mathcal{B}_{n-1}(z)+A\mathcal{B}_{n}(z)\,.
\end{array}
\right\}
$$
Substituting in ~(\ref{8}), and identifying with zero the
coefficients of the vectorial polynomials in the obtained
expression, we arrive to
\begin{equation}
\left.
\begin{array}{rcl}
\dot{B}_{n} & = & AC_{n}-C_{n-1}A+D_{n-1}B_{n}-B_{n}D_{n-1}\\
\dot{C}_n& = & D_nC_{n}-C_nD_{n-1}+B_{n+1}C_n-C_nB_{n}
\end{array}
\right\}\,n=1,2,\ldots  \label{B}
\end{equation}
Taking into account that, with the above notation, $D_n=\left(J_-\right)_{n+1,n+1}$, we see that ~(\ref{B}) is
equivalent to ~(\ref{Laxpair}) when we consider $J$ as a blocked matrix. \hfill $\square$

\section{Proof of Theorem~\ref{equivalencia} and Corollaries} \label{sec:4}

\subsection{Proof of Theorem~\ref{equivalencia}}

We start by establishing the equivalence between ~(\ref{bloques})
and ~(\ref{derivada bloques}). The key is the convergence in the
respective operator norm, for $|z|>\|J\|$, of the series given in
the right hand side of ~(\ref{resolvente}) and
~(\ref{resolventefinita}). Starting by ~(\ref{bloques}), to obtain
$\dot{\mathcal{R}}_J(z)$ it is sufficient to take derivatives in
~(\ref{resolventefinita}) and to substitute~$\dot{J}_{11}^n$~in
$$
\dot{\mathcal{R}}_J(z)=  \sum_{n\geq 0} \frac{\dot{J}^n_{11} }{z^{n+1}}\,,\quad |z| > \| J \|\,.
$$
Reciprocally, if ~(\ref{derivada bloques}) holds, substituting $\mathcal{R}_J(z)$ and $\dot{\mathcal{R}}_J(z)$ by
its Laurent expansion, and comparing their coefficients, we arrive to ~(\ref{bloques}).

The rest of the proof is to show the equivalence between (a), (b) and (d).

(a)$\Rightarrow$ (b) is Lemma~\ref{lema2}.

Now, we are going to prove (b) $\Rightarrow$ (d). We assume that ~(\ref{bloques}) is verified. Taking $n=1$ in
this expression we immediately deduce $\dot{a}_1=b_1$. Moreover, from ~(\ref{bloques}) we arrive to
~(\ref{derivada momentos}) in the same way that in the proof of (e)$\Rightarrow$ (f) in Theorem~\ref{teorema2}.
Then we are under the hypothesis of Theorem~\ref{teorema2} and, therefore, we have ~(\ref{derivada polinomios}).

Finally, ~(\ref{derivada polinomios}) $\Rightarrow$ ~(\ref{Laxpair}) was proved in Theorem~\ref{teorema2}. \hfill $\square$

\subsection{Proof of Corollary \ref{coro1}}

It is easy to see that $C_0$, $M(t,z)$ and $N(t)$ are, respectively, the solutions of the following Cauchy problems:
 \begin{multline*}
\left.
\begin{array}{r}
\dot{X}=-\left(J_-\right)_{11}X \\
X(0)=C_0(0)
\end{array}
\right\},\,
\left.
\begin{array}{r}
\dot{X}=-\exp(-zt)C_0^{-1}N(t) \\
X(0)=\left(C_0(0)\right)^{-1}\mathcal{R}_0(z)
\end{array}
\right\},\, \\
\text{and } \
\left.
\begin{array}{r}
\dot{X}=\left(J_{11}-\left(J_-\right)_{11}\right)X \\
X(0)=I_2
\end{array}
\right\} \, .
 \end{multline*}
Taking derivatives in the right hand side of ~(\ref{representacion}), and checking the initial condition, we can prove that $\exp(zt)C_0M(t,z)\left(N(t)\right)^{-1}$ is a solution of the following Cauchy problem,
\begin{equation}
\left.
\begin{array}{r}
\dot{X}=X\left(zI_2-J_{11}\right) -I_2+\left[X, \left(J_-\right)_{11}\right] \\
X(0)=\mathcal{R}_0(z)
\end{array}
\right\}\,.\label{44}
\end{equation}
 From ~\cite{Vilenkin}, we know that ~(\ref{44}) has a unique solution.
On the other hand, since Theorem~\ref{equivalencia}, $\mathcal{R}(z)$ is solution of ~(\ref{44}). Then, we arrive to ~(\ref{representacion}). \hfill $\square$

\subsection{Proof of Corollary \ref{corolario2}}

Since ~(\ref{momento0}), in the conditions of Corollary \ref{corolario2}, necessarily ~(\ref{Uexponencial}) implies
$$M=
\left[
\left(e^{zt} \, \mathcal{U}_0 \right)(\mathcal{P}_0)
\right]^{-1}\,.
$$
On the other hand, for proving  that
$\{a_n\,,\,b_n\,,\,c_n\}\,,\,n\in
\mathbb{N}\,,$ is a solution of~(\ref{Kostant}), it is sufficient to show that ~(\ref{derivada momentos}) holds. Let a fixed $k\in\{0,1,\ldots, \}$ be. Taking into account
$$
\frac{d}{dt}\left(e^{zt} \, \mathcal{U}_0\right)(z^k\mathcal{P}_0)=\left(e^{zt} \, \mathcal{U}_0\right)(z^{k+1}\mathcal{P}_0)
$$
and
$$
\frac{d}{dt}M=-M\mathcal{U}(z\mathcal{P}_0)\,,
$$
and taking derivatives in
$$
\mathcal{U}(z^k\mathcal{P}_0)=
\left(e^{zt} \, \mathcal{U}_0\right)(z^k\mathcal{P}_0)M \, ,
$$
see ~(\ref{Uexponencial}), we arrive to ~(\ref{derivada momentos}). \hfill $\square$

%\section*{Acknowledgments}
%This research was carried out while on a visit of the first author at University of Coimbra. The work of this author was supported in part by Direcci\'{o}n General de Investigaci\'{o}n, Ministerio de Educaci\'{o}n y Ciencia, under grant MTM2006-13000-C03-02. The work of the second author was supported by CMUC/FCT. The third author would
%like to thank UI Matem\'atica e Aplica\c c\~oes from University of Aveiro.

\end{document}